\documentclass{amsart}

\newtheorem{theorem}{Theorem}[section]
\newtheorem{lemma}[theorem]{Lemma}
\newtheorem{proposition}[theorem]{Proposition}
\newtheorem{corollary}[theorem]{Corollary}

\theoremstyle{definition}
\newtheorem{definition}[theorem]{Definition}
\newtheorem{example}[theorem]{Example}
\newtheorem{ack}{Acknowledgements}

\theoremstyle{remark}
\newtheorem{remark}[theorem]{Remark}

\newcommand{\A}{{\mathcal A}}

\newcommand{\bnbc}{\beta{\rm {\bf nbc}}}

\newcommand{\C}{\mathbb C}
\newcommand{\codim}{\mbox{\rm codim}}
\newcommand{\ints}{\mathbb Z}

\newcommand{\calL}{{\mathcal L}}
\newcommand{\la}{\lambda}
\newcommand{\calF}{{\mathcal F}}
\newcommand{\nbc}{{\rm {\bf nbc}}}
\newcommand{\om}{\omega}

\newcommand{\stR}[1]{\stackrel{#1}{\longrightarrow}}

\newcommand{\w}{\omega}
\newcommand{\we}{{\wedge}}

\begin{document}

\title[$\beta${\nbc}-bases for cohomology of local systems]{
$\beta${\bf nbc}-bases for cohomology of local systems on
hyperplane complements}
\author{Michael Falk}
\thanks{first author partially supported by NAU Organized Research
Grant}
\address{Department of Mathematics \\ Northern Arizona
University \\ Flagstaff, AZ 86011}
\email{mjf@odin.math.nau.edu}
\author{Hiroaki Terao}
\address{Department of Mathematics\\
University of Wisconsin - Madison \\ Madison, WI 53704}

\email{terao@math.wisc.edu}

\dedicatory{In memory of Michitake Kita}


\begin{abstract} {\sl \small We study cohomology with coefficients in a
rank one local system on the complement of an arrangement of
hyperplanes
$\A$. The cohomology plays an important role for the theory of
generalized hypergeometric functions. We combine several known results
to construct explicit bases of logarithmic forms for the only
non-vanishing cohomology group, under some nonresonance conditions on
the local system, for any arrangement $\A$. The bases are determined by
a linear ordering of the hyperplanes, and are indexed by certain
``no-broken-circuits" bases of $\A$. The basic forms depend on the
local
system, but any two bases constructed in  this way are related by a
matrix of integer constants which depend only on the linear orders and
not on the local system. In certain special cases we show the existence
of bases of monomial logarithmic forms. }
\end{abstract}

\maketitle

\begin{section}{Introduction}

Let $V$ be the $\ell$-dimensional affine space over $\C$, and let $\A$
be an arrangement of hyperplanes in $V$. We  fix a linear order on $\A$
and write $\A=\{H_1, \ldots, H_n\}$. Let $M=M(\A)=V - \bigcup_{i=1}^n
H_i$ denote the complement of $\A$. Let $L=L(\A)$ be the intersection
poset of $\A$. By definition $L$ is the set of nonempty intersections
of
hyperplanes in $\A$ ordered by reverse inclusion. By convention $L$
includes $V$ as its unique minimal element. Then $L$ is a ranked poset,
with $r(X)=\codim(X)$, and all maximal elements have the same rank
\cite[Lemma 2.4]{OrT1}. The rank of $\A$ is then the rank of any
maximal
element of $L$, denoted by $r$.

Let $\alpha_i$ be a polynomial function of degree one on $V$ which
vanishes exactly on $H_i$, for $1 \leq i \leq n$. Given a complex
$n$-vector $\la=(\la_{1}, \ldots, \la_n)$ consider the multi-valued
holomorphic function $U_\la$ defined on $M$ by $U_\la=\alpha_1^{\la_1}
\cdots \alpha_n^{\la_n}$. Let $\calL_\la$ be the rank one local system
on $M$ whose local sections are constant multiples of branches of $1 /
{U_\la}$.

The cohomology $H^r(M,\calL_\la)$ is important for the theory of
integrals of $U_\la$, the Aomoto-Gelfand theory of generalized
hypergeometric functions, which itself has many important applications
in diverse areas of mathematics. Consult \cite{Gel1,AoK1,Var2} for
references.
In this context the cohomology $H^q(M,\calL_\la)$ was the main subject
of \cite{Aom2,Koh1,ESV1}, among others.

It is known that $H^q(M,\calL_\la)$ vanishes for $q \not = r$ under
certain nonresonance conditions on $\la$. The main goal of this note is
to find an explicit basis for $H^r(M,\calL_\la)$ for arbitrary $\A$ and
generic $\la$. This basis depends on the choice of linear order on $\A$
and is, in some sense, the analogue of the {\bf nbc} (``no broken
circuits") basis for the ordinary cohomology $H^r(M,\C)$
\cite{B1,OrT1}.
The basis is parametrized by the $\bnbc$-bases of $\A$, as defined by
Ziegler in \cite{Zie1}. Our main result is obtained by combining the
results of \cite{ESV1}, \cite{Yuz1}, and \cite{Zie1}.

When $\A$ is a
general position arrangement, with real defining
forms, our basis coincides
with the basis which was first found by Aomoto \cite[p.292]{Aom1}.
If $\A$ is a normal arrangement
or in general position to infinity (both defined in section 4), again
with real defining forms, then our
basic forms coincide with forms constructed by Varchenko
\cite[6.2]{Var1} -- see Remark \ref{Varch}.
The construction in this paper
applies in complete generality, without the restriction to real
defining
forms. The resulting basis is determined by the underlying
combinatorial
structure (the affine matroid), whereas the bases constructed in
\cite{Aom1,Var1} depend on the geometry of the associated real
arrangement.
In independent work, A.~Douai \cite{Dou1} has produced a basis
for $H^r(\A)$ for general
$\A$ and generic (i.e. transcendental) $\la$. The methods and the
resulting basis are different from ours.
In addition, construction of a basis for $H^r(\A)$ forms a part of the
thesis of H.~Kanarek \cite{Kan1}.

For real arrangements Varchenko \cite{Var1} conjectures a formula for
the determinant of a
matrix of hypergeometric integrals $\int_{\Delta_i}
U_\la\omega_{\Delta_j}$
whose rows are indexed by
the bounded chambers $\Delta$ of the real form of $\A$, and
whose columns are indexed by certain hypergeometric forms
$\omega_\Delta \in H^r(M,\calL_\la)$.
There is a natural one-to-one correspondence between bounded chambers
and $\bnbc$-bases of $\A$. So our result provides a variation
of Varchenko's conjecture in which the forms $\omega_\Delta$ ar
replaced
by our basic forms $\Xi_B$. These forms are easier to describe and
satisfy a simple recursion, perhaps facilitating a proof of Varchenko's
conjecture. This alternate formulation is adopted in \cite{Dou1}.
In the cases where the determinant formula has been proven, the forms
$\omega_\Delta$ coincide with the $\Xi_B$, but not in
general.

In the definition of generalized hypergeometric functions, the
arrangement $\A$ is considered to be an independent variable varying in
over a parameter space $\Gamma$ of arrangements with constant
underlying
matroid. The cohomology groups $H^r(M,\calL_\la)$ comprise a vector
bundle over $\Gamma$. Because of the combinatorial nature of our
construction, we actually get a complete system of independent global
sections, i.e., a global frame,
of this vector bundle. This may allow one to give more explicit
formulas for the associated generalized hypergeometric functions.

In the case of affine supersolvable arrangements, defined in Remark
\ref{affss}, it is particularly
easy to write down the $\bnbc$-bases of $\A$. This case includes the
discriminantal arrangements of \cite{ScV1}, also called Selberg-type
arrangements in \cite{AoK1}. These arrangements are parametrized by the
configuration space consisting of sets of $n$ pairwise distinct point
in
$\C^n$; thus the arrangements do not in general have real defining
forms.
Schechtmann and Varchenko use this setup to show that solutions of the
Knizhnik-Zamolodchikov equations of conformal field theory are given by
generalized hypergeometric integrals \cite{ScV1}. The $\bnbc$-bases of
global sections may then be used to construct explicit solutions of the
KZ equations.
Along the same lines, such a basis may be used to study the monodromy
action on the twisted cohomology, in particular for calculation of the
monodromy of the KZ equation. This is one focus of the dissertation
\cite{Kan1}.

Unless otherwise noted, we adopt the notation of \cite{OrT1}, which
should also be consulted for definitions and background material on
arrangements.

In section 2 we review the relevant combinatorial constructions as
developed in \cite{Zie1}. These notions have been reformulated in more
geometric language, and all references to the associated central
arrangement $c\A$ have been eliminated. We define the broken circuit
complex $BC$ of $\A$, whose simplices are $\nbc$ sets of $\A$. The
lexicographic order on the facets of $BC$ is a shelling order, and
leads
immediately to the notion of $\bnbc$-base. We observe that the
$\bnbc$-bases yield a natural basis for the cohomology
$\tilde{H}^{r-1}(BC)$.
The $\bnbc$-bases of $\A$ can also be constructed inductively using
deletion-restriction \cite[Thm 1.5]{Zie1}, and we describe this
process.
Example \ref{example} illustrates these ideas.

Let $F$ be the Folkman complex of $\A$, defined as the simplicial
complex of chains in $L - \{V\}$. There is a natural homotopy
equivalence of $BC$ with the Folkman complex $F$, which we use to
construct a basis for $\tilde{H}^{r-1}(F)$.

In section 3 we review and study the results of Yuzvinsky \cite{Yuz1}.
Define logarithmic differential forms $\w_i=
\frac{d\alpha_i}{\alpha_i}$
on $V$. The Orlik-Solomon algebra $A^\cdot$ of $L$ is the
finite dimensional graded $\C$-algebra generated by the $\w_i$ under
exterior product. By convention $A^0=\C$. Define a logarithmic 1-form
$\w_\la$ by $$\w_\la=d(\log(U_\la))=\sum_{i=1}^n \la_i\w_i.$$ Also
define a map $\w_\la \we : A^p \to A^{p+1}$ by $w_\la \we
(\eta):=\w_\la
\we \eta.$ Then $(\w_\la \we)^2=0$ so $(A^\cdot,\w_\la \we)$ is a
complex. In \cite{Yuz1} Yuzvinsky showed that there exists a natural
isomorphism $$\phi: H^r(A^\cdot,\w_\la \we) \to \tilde{H}^{r-1}(F),$$
using sheaf theory on posets, under certain genericity conditions on
$\la$. This shows in particular that $H^q(A^\cdot,\w_\la \we)=0$ for
$q<r$.

We explicitly describe the inverse map $\phi^{-1}$, for $q=r$
(Proposition \ref{phiinverse}).
The map $\phi^{-1}$ turns out to coincide up to sign with
the one induced by the map $S^{r} $ introduced and studied by
Schechtman and Varchenko in \cite[3.2]{ScV1}.
The image of the $\bnbc$-basis for
$\tilde{H}^{r-1}(F)$ under $\phi^{-1}$ yields a basis for
$H^r(A^\cdot,\w_\la \we)$. On the other hand, Esnault, Schechtman, and
Viehweg proved in \cite{ESV1} that the de Rham map
$$H^r(A^\cdot,\w_\la \we) \to H^r(M,\calL_\la)$$ is an isomorphism
under
appropriate nonresonance conditions on $\la$. A more general version of
this result, with fewer restrictions on $\la$, appears in \cite{STV1}.
Thus we obtain a basis of logarithmic forms for $H^r(M,\calL_\la)$. The
explicit description of this basis is in Theorem \ref{main}.

More important than the result itself is the method of proof.
Yuzvinsky's isomorphism allows questions about the local system
cohomology  to be carried over to the broken circuit complex, where
powerful combinatorial tools apply. In particular, the parameter $\la$
does not appear in $BC$. We can conclude that two bases constructed as
above from different linear orderings of $\A$ are connected by a
transition matrix which is independent of $\la$, and in fact has
integral entries.

One might ask if the monomials $\w_B=\prod_{k=1}^r \w_{i_k}$, for
$B=(H_{i_1}, \ldots, H_{i_r})$, $i_{1}<\cdots < i_{r}$, themselves form
a basis, as $B$ ranges over the $\bnbc$-bases of $\A$. In section 4 we
show that this is the case, under the same nonresonance conditions on
$\la$, when the linear ordering on $\A$ is {\em unmixed} (Def.
\ref{unmixeddef}). If $\A$ is in general position, or, more generally,
if the hyperplane at infinity is generic relative to $\A$, then every
linear order is unmixed.
This is also the case when
$\bigcup\A$ is a normal crossing divisor.
There are some other special linear orders on
arrangements of rank two for which these
monomials form a basis.  But
without imposing some additional unnatural genericity requirements on
$\la$, the monomials above will not form a basis for arbitrary orders,
even in rank two.

\end{section}

\begin{section}{Broken circuits and $\beta$-systems}

In this section we establish some notation and develop the
combinatorial
tools needed for the proof of the main theorem in section 3. All this
material is adapted from the references \cite{BjZ1} and \cite{Zie1}.

Let $V$ be the $\ell$-dimensional affine space over $\C$. Let $\A$ be
an
arrangement of hyperplanes in $V$. Fix a linear order on $\A$ by
labelling the hyperplanes $H_1, \ldots, H_n$. We will sometimes denote
this linear order by $<$ when no confusion will result. The {\em
intersection poset\/} $L$ of $\A$ consists of the nonempty affine
subspaces $X= \bigcap_{i \in I} H_i$ for $I \subseteq \{1, \ldots n\}$,
partially ordered by reverse inclusion. We will occasionally write
$H_I$
for $\bigcap_{i \in I}H_i$. The ambient space $V$, corresponding to
$I=\emptyset$, is the unique minimal element of $L$. The maximal
elements of $L$ all have the same codimension \cite[Lemma 2.4]{OrT1},
which we denote by $r$. For $X \in L$ we set $\A_X=\{H \in \A \ | \ X
\subseteq H\}$.

A subset $\{H_i \ | \ i \in I\}$ of $\A$ is {\em
dependent\/} if $\bigcap_{i \in I} H_i \not = \emptyset$ and
$\codim(\bigcap_{i \in I} H_i)<|I|$. A subset of $\A$ which has
nonempty
intersection and is not dependent is called {\em independent}. Maximal
independent sets are called {\em bases}; every base has cardinality
$r$.

An inclusion-minimal dependent set is called a {\em circuit}. A {\em
broken circuit\/} is a subset $S$ of $\A$ for which there exists $H <
\min(S)$ such that $\{H\} \cup S$ is a circuit. The collection of
subsets of $\A$ having nonempty intersection and containing no broken
circuits is a simplicial complex we call the {\em broken circuit
complex\/} of $\A$, denoted by $BC$. This is a pure
$(r-1)$-dimensional complex consisting of independent sets. Simplices
of
$BC$ are called \nbc\ sets, and facets (maximal simplices) of $BC$ are
bases of $\A$ called \nbc-bases. We say an ordered base $(H_{i_1},
\ldots, H_{i_r})$ is {\em standard\/} if $i_1<\ldots <i_r$. We denote
by
\nbc\ the set of standard ordered \nbc-bases of $\A$.

The {\em Folkman complex\/} of $L$ is the simplicial complex of
linearly
ordered subsets $(X_0 > \ldots > X_p)$ of $L - \{V\}$. It is also a
pure
$(r-1)$-dimensional simplicial complex, which we denote by $F$.

The purpose of this section is to expose a natural basis for
$\tilde{H}^{r-1}(F)$, given the linear ordering on $\A$. First we
describe the
structure of $BC$, following \cite{Zie1}. A {\em shelling order\/} for
a
simplicial complex $K$ is a linear order $\prec$ on the set $\mathcal
F$ of
facets of $K$ such that $\sigma \cap (\bigcup_{\tau \prec \sigma}
\tau)$
is a nonempty union of facets of $\sigma$, for each $\sigma \in
{\mathcal F}
- \min({\mathcal F})$.
Suppose $\prec$ is a shelling order for $K$. A facet $\sigma$ is called
a {\em homology facet\/} if $\sigma \cap (\bigcup_{\tau \prec \sigma}
\tau)=\partial \sigma$. It is easy to see that the union of the
non-homology facets of $K$ is a contractible subcomplex, and thus $K$
has the homotopy type of a bouquet of spheres, one for each homology
facet.

\begin{theorem}[{\rm \cite[Thm.7.4.3]{B2}}]{The lexicographic ordering
of \nbc\ is a
shelling order for $BC$.} \label{shelling}
\end{theorem}

If $B \in \nbc$ and $H \in B$ then the facet $B \-- \{H\}$ of $B$ lies
in $B \cap (\bigcup_{B'\prec B} B')$ if and only if $B':= (B \-- \{H\})
\cup
\{H'\}$ is an \nbc-base for some $H'<H.$
By \cite[Lemma 1.2]{Zie1} it is enough that $B'$ is a base for some
$H'<H$.
This determines the homology facets
in the shelling of Theorem \ref{shelling}.

\begin{definition}{A base $B$ is called a {\em $\bnbc$-base} if $B$
is an \nbc-base and, for every  $H\in B$ there exists $H'<H$ in $\A$
such that $ (B - \{H\}) \cup \{H'\} $ is a base. }
\end{definition}

\begin{theorem}[{\rm Ziegler}]{The
homology facets for the lexicographic
shelling of $BC$ are precisely the $\bnbc$-bases of $\A$.}
\label{facets}
\end{theorem}

The set of $\bnbc$-bases of $\A$ is called a {\em $beta$-system\/} for
$\A$. The set of standard ordered $\bnbc$-bases of $\A$ will be denoted
here by $\bnbc$. The notation and terminology comes from the fact that
the cardinality of $\bnbc$ is equal to Crapo's beta invariant of the
matroid of $c\A$. In case $r=\ell$ and $\A$ is complexified real this
is
precisely the number of bounded chambers of the real form of $\A$.

There is also an inductive (``deletion-restriction") definition of
$\bnbc$ provided by \cite[Thm. 1.5]{Zie1}. We say $H \in \A$ is a {\em
separator\/} if the rank of $\A - \{H\}$ is less than $r$. Let
$(\A,\A',\A'')$ denote a triple of arrangements relative to the
hyperplane $H_n$, as defined in \cite{OrT1}. For $H'' \in \A''$ let
$\nu(H'')$ denote the smallest hyperplane of $\A$ containing $H''$.
(Note
that $\nu$ is injective.) We may apply $\nu$ to $(r-1)$-tuples of
hyperplanes in the obvious way. Order $\A'$ as a subset of $\A$. Order
$\A''$ by $H''<K''$ if and only if $\nu(H'')<\nu(K'')$.

\begin{theorem}[{\rm Ziegler}]{If $H_n$ is a separator then
$\bnbc(\A)=\emptyset$. Otherwise $$\bnbc(\A)=\bnbc(\A') \cup
\{(\nu(B''),
H_n) \ | \ B'' \in \bnbc(\A'')\}.$$}
\label{induct}
\end{theorem}

\begin{example} {If $\A$ is an arrangement of rank one, then
$\nbc=\A$ and $\bnbc=\A \-- \min(\A)$.
Let $\A$ be the
arrangement of rank two with defining equation
$$Q=(x+1)(x-1)(y+1)(y-1)(x-y),$$ with the linear order defined by the
order of the factors in $Q$. Then
$$\nbc=\{(H_1,H_3),(H_1,H_4),(H_1,H_5),(H_2,H_3),(H_2,H_4),(H_2,H_5)\}$$
and $$\bnbc=\{(H_2,H_4),(H_2,H_5)\}.$$ Here $\A'=\{H_1,H_2,H_3,H_4\}$
with $\bnbc(\A')=\{(H_3,H_4)\}$, and $\A''=\{H_1\cap H_5,H_2\cap H_5\}$
so that $\bnbc(\A'')=\{H_{245}\}$ and $\nu(H_2\cap
H_5)=\min\{H_2,H_4,H_5\} =H_2$.}
\label{example}
\end{example}

\begin{remark}{The previous example is a special case of a more
general phenomenon. Let us say $\A$ is {\em supersolvable\/} if there
is a
sequence $\A_1 \subset \ldots \subset \A_r=\A$ such that the rank of
$\A_p$ is equal to $p$ and, for distinct $H,H' \in \A_p$ with
$H \cap H'\not = \emptyset$, there exists $H''\in \A_i$ with $i<p$ and
$H'' \supset H\cap H'.$ Examples of supersolvable affine arrangements
include the discriminantal arrangements of \cite{ScV1}.

Suppose $\A$ is supersolvable. Set
${\mathcal B}_p=\A_p \-- \A_{p-1}$. Order $\A$ so that the elements of
${\mathcal B}_{p-1}$ precede those of ${\mathcal B}_p$ for $1<p\leq
r$.
It is proved in
\cite{BjZ1} that $B=(H_{i_1}, \ldots, H_{i_r}) \in \nbc$ if and
only if $H_{i_p}\in
{\mathcal B}_p$ for all $p$. One immediately concludes that
$B\in \bnbc$ if and only if
$H_{i_p}\in {\mathcal B}_p \-- \min {\mathcal B}_p$ for all $p$.}
\label{affss}
\end{remark}

For $B \in \nbc$ let $B^* \in C^{r-1}(BC)$ denote the cochain dual to
$B$. Thus $B^*$ is determined by the formula $$\langle B^*,
B'\rangle=\left\{ \begin{array}{ll} 1 & \mbox{if  $B'=B$} \\ 0 &
\mbox{otherwise}\\ \end{array}\right. $$ for $B' \in \nbc$.

The next proposition follows immediately from Theorem \ref{facets},
using the properties of shellable complexes outlined earlier.

\begin{proposition}{The set $\{[B^*] \ | \ B \in \bnbc\}$ forms a basis
for
$\tilde{H}^{r-1}(BC)$.}
\label{bcbasis}
\end{proposition}

The next step is to carry the basis of Proposition \ref{bcbasis} over
to
$F$. The connection is provided by the following result.

\begin{theorem}[{\rm Bj\"orner, Ziegler}]{The vertex map $X \mapsto
\min(\A_X)$ induces a simplicial map $\pi: F \longrightarrow BC$, and
$\pi$ is a homotopy equivalence.}
\label{htpyeq}
\end{theorem}

\begin{proof} This is a special case of Theorem 3.12 of \cite{BjZ1}.
\end{proof}

\noindent
It follows that $F$ has the homotopy type of a bouquet of
$(r-1)$-spheres \cite{WaW}.

Combining Theorem \ref{htpyeq} and Proposition \ref{bcbasis} we have
that
$\{\pi^*([B^*]) \ | \ B \in \bnbc\}$ forms a basis for
$\tilde{H}^{r-1}(F)$. We
need a slightly more precise result. For an $(r-1)$-simplex $\xi \in
F^{(r-1)}$ let $\xi^* \in C^{r-1}(F)$ be the cochain dual to $\xi$. If
$B=(H_{i_1}, \ldots, H_{i_r})$ is an ordered base, set $$\xi(B)=(X_1 >
\ldots > X_r)$$ where $X_p=\bigcap_{k=p}^r H_{i_k}$, for $1 \leq p \leq
r.$

\begin{theorem}{The set $\{[\xi(B)^*] \ | \ B \in \bnbc\}$ forms a
basis
for $\tilde{H}^{r-1}(F)$.}
\label{fbasis}
\end{theorem}

The proof will follow from Corollary \ref{pi*} below. \begin{lemma}{Let
$B=(H_{i_1}, \ldots , H_{i_r})$ be a standard ordered base, with
$\xi(B)=(X_1>\ldots > X_r)$. Then $B \in \nbc$ if and only if
$H_{i_k}=\min \A_{X_k}$ for all $1 \leq k \leq r$.}
\label{min}
\end{lemma}

\begin{proof} The base $B$ contains a broken circuit with smallest
element $H_{i_k}$ if and only if $H_{i_k} \not = \min \A_{X_k}.$
\end{proof}

\begin{lemma} \label{inject} {Let $B=(H_{i_1}, \ldots ,H_{i_r}) \in
\nbc$.
\begin{enumerate}
\item[{\rm (i)}] Let $\xi \in F^{(r-1)}$. Then $\pi(\xi)=B$ if and only
if
$\xi=\xi(B)$.
\item[{\rm (ii)}] Let $B'=(H'_{i_1},\ldots ,H'_{i_r})$ be any ordered
\nbc-base.
Then $\xi(B)=\xi(B')$ if and only if $B=B'$.
\end{enumerate} }
\end{lemma}

\begin{proof} (i) Sufficiency follows immediately from Lemma
\ref{min}.
Assume $\pi(\xi)=B$ and write $\xi=(X_1 > \ldots > X_r)$. Then
$\A_{X_p}\supset \A_{X_{p+1}}$ for $1\leq p<r$. Then $H_{i_p}\cap
\ldots
\cap H_{i_r}\supseteq X_p$ for all $1\leq p\leq r$, and it follows from
a dimension count that $\xi=\xi(B)$.

(ii) Assume $\xi(B)=\xi(B')$. Suppose $B \not =B'$ and let $p$ be
maximal with $H_{i_p} \not = H'_{i_p}$. Since $H_{i_p}\cap \ldots
H_{i_r}=X_p=H'_{i_p}\cap \ldots H'_{i_r}$ by assumption, this implies
that $\{H_{i_p}, H'_{i_p},\ldots,H'_{i_r}\}$ is dependent. It then
follows from Lemma \ref{min} that $\{H'_{i_p},\ldots ,H'_{i_r}\}$
contains a broken circuit, and this is a contradiction.
\end{proof}

\noindent
Let $\pi^\sharp: C^{r-1}(BC) \to C^{r-1}(F)$ be the map on the
cochain level induced by $\pi$.

\begin{corollary}{For any $B \in \nbc$, $\pi^\sharp(B^*)=\xi(B)^*$.}
\label{pi*}
\end{corollary}

\begin{proof} By definition $\langle \pi^\sharp(B^*),\xi
\rangle=\langle
B^*,\pi(\xi) \rangle$. By Lemma \ref{inject} (i), $\langle B^*,\pi(\xi)
\rangle= \langle \xi(B)^*,\xi\rangle,$ proving the equality of the
lemma.
\end{proof}

\noindent Now Theorem \ref{fbasis} follows immediately from Corollary
\ref{pi*} together with Proposition \ref{bcbasis} and Theorem
\ref{htpyeq}.

\end{section}

\begin{section}{$\beta${\bf nbc}-bases for $H^{r}(M, \calL_{\la})$} We
first introduce a topology on the intersection poset $L=L(\A)$ as
follows: $O \subseteq L$ is {\em open} if and only if $O$ is a
lower ideal, that is, if $X\in O$ and $Y\leq X$, then $Y\in O$. For
$X\in
L$ let $L_{X} =\{Y\in L | Y \leq X\}$. Then $L_{X} $ is an open set.

Let $0\leq p\leq r$.  Define a sheaf $\calF^{p}$ on $L$  whose stalk
$\calF^{p}(X) $ at $X\in L$  is defined to be equal to $A^{p}(\A_{X}
)$,
the degree $p$ part of the Orlik--Solomon algebra $A^{\cdot}(\A_{X})$
of
the arrangement $\A_{X}.$ If $Z\in L$ and $\codim Z=p$, put $A_{Z} =
A_{p}(\A_{Z}).$ Then $A_{p}(\A_{X} ) = \oplus A_{Z},$ summing over all
$Z\in L$ of codimension $p$ such that $Z\leq X$ \cite[3.73]{OrT1}.
Thus
if $Y\leq X$, then $\calF^{p}(Y) $  is a direct summand of
$\calF^{p}(X)$ and there is a projection $$ \rho_{X, Y} : \calF^{p}(X)
\longrightarrow \calF^{p}(Y). $$ These projections define a sheaf
$\calF^{p}$ on $L$, with $\Gamma (L_{X}, \calF^{p}) = A_{p}(\A_{X})$
\cite[3.2]{Yuz1}.

The results of this section require that $\la$ satisfy certain
genericity conditions. These are given below in Remark \ref{conv}. At
this point we must assume $\la_i\not =0$ for all $i$. For each $X\in
L$,
let $$ \om_{\la}(X) = \sum_{H_i\in \A_X} \la_{i}\om_{i}.$$ Define a
sheaf homomorphism \[ d_{\la} : \calF^{p} \longrightarrow \calF^{p+1}
\]
by defining \[ d_{\la}(X) : \calF^{p}(X) \longrightarrow
\calF^{p+1}(X),
\] as the left multiplication by $\om_{\la}(X)$, for each $X\in L$. The
sheaf $\calF^{0} $ is the constant sheaf with $\calF^{0}(X) = \C$ for
any $X\in L$. Then $\ker (d_{\la}) \subset \calF^{0}$ is the skyscraper
sheaf $\C_{V}$ whose only nonzero stalk is at $V$ and it equals $\C$.
Thus we have a complex of sheaves \[ 0 \longrightarrow \C_{V}
\longrightarrow \calF^{0} \stR{d_{\la}} \calF^{1} \stR{d_{\la}} \cdots
\stR{d_{\la}} \calF^{r} \longrightarrow 0. \] Yuzvinsky \cite{Yuz1}
showed that this complex gives a flabby resolution of $\C_{V}$ under
certain genericity conditions on $\la$. In order to formulate these
non-resonance conditions we need the following definition.

\begin{definition} \label{dense}{An element $X\in L - \{V\}$ is {\em
dense} if $\A_{X}$ is not decomposable \cite{STV1}, that is, if
$\A_{X}$
is not a product of two nonempty arrangements.}
\end{definition}

\noindent
In \cite{Yuz1} the condition $\bar{\chi}(A(\A_{X})) \neq 0$ is
considered.  This condition is actually equivalent to the denseness of
$X$ by \cite{Cra}. (Also see \cite{STV1}.)

\begin{remark} \label{conv} {For $X \in L$ set $\la(X)=\sum_{H_i
\in \A_X} \la_i$. For the remainder of this paper we assume $\la(X)\not
=0$ for all dense $X \in L \-- \{V\}$.  Note that this implies $\la_{i}
\neq 0$
for all $i$.}
\end{remark}

\begin{theorem}[{\rm Yuzvinsky}]{The complex of sheaves \[ 0
\longrightarrow \C_{V} \longrightarrow \calF^{0} \stR{d_{\la}}
\calF^{1}
\stR{d_{\la}} \cdots \stR{d_{\la}} \calF^{r} \longrightarrow 0 \] is a
flabby resolution of the skyscraper sheaf $\C_{V}$.}
\label{Yuzflabby}
\end{theorem}

Note $\Gamma(L, \calF^{p}) = A^{p}(\A)$.  By the standard sheaf
cohomology theory \cite{God1}, we have an isomorphism \[ \phi :
H^{p}(A^{\cdot},
\om_{\la}\wedge) \tilde{\longrightarrow} H^{p}(L, \C_{V}). \]

Let $\calF$ be any sheaf on $L$.  Applying the canonical simplicial
resolution from \cite[6.4]{God1}, we can view $H^{p}(L, \calF) $ as the
$p$th cohomology of the cochain complex $C^{\cdot}(L, \calF) :=
\oplus_{p=0}^{r} C^{p}(L, \calF)$, where $C^{p}(L, \calF)$ is the set
of
all functions on the $(p+1)$-- tuples $(X_{0} > X_{1} > \cdots >
X_{p})$
from $L$  with values in $\calF(X_{0})$. The differential of the
complex
is given by
\begin{eqnarray*} (df)(X_{0} > X_{1} > \cdots > X_{p})
&=& (-1)^{p} \rho_{X_{p-1}, X_{p}} f(X_{0} > \cdots > X_{p-1}) \\&&+
\sum_{j=0}^{p-1}  (-1)^{j}  f(X_{0} > \cdots > \widehat{X_{j}} > \cdots
> X_{p} ). \end{eqnarray*}
We often simply write $C^{p}(\calF)$ instead
of $C^{p}(L, \calF)$.

Let $F = F(\A)$ be the Folkman complex of $\A$, defined in section 2.
Let $C^{p-1}(F, \C)$ be group of simplicial $(p-1)$--cochains on $F$
for
$1\leq p\leq r$.  Define $C^{-1}(F, \C) = \C$ and consider the
augmented
cochain complex $C^{\cdot} (F, \C)$ with the usual coboundary maps.
Consider the cochain maps \[ \gamma^{p} : C^{p-1}(F, \C)
\longrightarrow
C^{p}(L, \C_{V}) ~~~~(0\leq p\leq r) \] given by
\[ (\gamma^{p}(f)) (X_{0}
> X_{1} > \cdots > X_{p}) = \left\{ \begin{array}{ll} f(X_{0} > \cdots
>
X_{p-1}) & \mbox{\rm if~} X_{p} = V,\\ 0 & \mbox{\rm otherwise}.
\end{array} \right. \]
for $f\in C^{p-1}(F, \C)$.
Then $\gamma =
\{\gamma^{p} \}$ defines an isomorphism between the two cochain
complexes. Using $\gamma$, we identify the sheaf cohomology $H^{p}(L,
\C_{V})$ with the reduced simplicial cohomology $\tilde{H}^{p-1}(F,
\C)$. Then we have the following corollary \cite[4.1]{Yuz1}.

\begin{corollary}[{\rm Yuzvinsky}] \label{Yuzcor}
\begin{enumerate}
\item $H^{p}(A^{\cdot}, \om_{\la}\wedge) \simeq \tilde{H}^{p-1}(F, \C)
=
0$ unless $p=r$, and
\item $H^{r}(A^{\cdot}, \om_{\la}\wedge) \simeq
\tilde{H}^{r-1}(F, \C).$
\end{enumerate}
\end{corollary}

\noindent
We denote the isomorphism \ref{Yuzcor} (ii) by
\[ \phi : H^{r}(A^{\cdot},
\om_{\la}\wedge) \tilde{\longrightarrow} \tilde{H}^{r-1}(F, \C). \]
We will describe the inverse isomorphism
$\phi^{-1}$ explicitly.

Let $\xi = (X_{1} > \cdots > X_{r}) \in F^{(r-1)} $. Let
$\xi^{*} \in C^{r-1}(F, \C)$ be the cochain dual to $\xi.$
Define a
linear map \[ \upsilon : C^{r-1}(F, \C) \longrightarrow A^{r} \] by \[
\upsilon (\xi^{*}) := \om_{\la}(X_{1})\cdots \om_{\la}(X_{r}). \]

\begin{remark}
\label{}
{
The form $\om_{\la}(X_{1})\cdots \om_{\la}(X_{r})$ is called
the flag form associated to the flag $\xi$, which was introduced and
studied by
Varchenko in \cite[6.1]{Var1}.  (See also \cite[3.2.3]{ScV1} and
\cite[10.2.11]{Var2}.)
}
\end{remark}

\begin{proposition}
\label{phiinverse}
Suppose $\la(X) \neq 0$ for all
dense $X\in L - \{V\}$. Then \[ \phi^{-1}:
\tilde{H}^{r-1} (F, \C) \longrightarrow H^{r}(A^{\cdot},
\om_{\la}\wedge) \] is induced by
$(-1)^{\frac{r(r+1)}{2}} \upsilon$.
\end{proposition}

\begin{proof}
The isomorphism $\phi$ is constructed out of a sequence of
connecting homomorphisms.
We have to trace back through this sequence.
Let $\xi=(X_{1} > \cdots > X_{r})\in F^{(r-1)}$. Recall
$\xi^{*} \in C^{r-1}(F, \C)$. Let $\beta_{0}$ be the corresponding
element
in $ C^{r}(\C_{V}) \subset C^{r} (\calF^{0})$. Then \[ \beta_{0} (Y_{1}
>\cdots > Y_{r} > V) = \left\{ \begin{array}{ll} 1 & \mbox{\rm if~}
(Y_{1} > \cdots > Y_{r} ) = (X_{1} > \cdots > X_{r})
 \\ 0 & \mbox{\rm otherwise}. \end{array}
\right. \] Define $\beta_{1} \in C^{r-1} (\calF^{0})$ by \[ \beta_{1}
(Y_{1} > \cdots > Y_{r} ) = \left\{ \begin{array}{ll} (-1)^{r} &
\mbox{\rm if~} (Y_{1} > \cdots > Y_{r} ) =
(X_{1} > \cdots > X_{r}) \\ 0 & \mbox{\rm
otherwise}. \end{array} \right. \] Then obviously $d\beta_{1} =
\beta_{0}.$ Let $\beta_{2} = d_{\la} \beta_{1} \in C^{r-1}
(\calF^{1})$.
Then \[ \beta_{2} (Y_{1} > \cdots > Y_{r} ) = \left\{ \begin{array}{ll}
(-1)^{r} \om_{\la}(X_{r}) & \mbox{\rm if~} (Y_{1} > \cdots > Y_{r} ) =
(X_{1} > \cdots > X_{r})
\\ 0 & \mbox{\rm otherwise}. \end{array} \right. \] Repeating this
construction, we obtain $\beta_{3} \in C^{r-2}(\calF^{1})$, $\beta_{4}
\in C^{r-2}(\calF^{2}),$ and so on, until finally we reach $\beta_{2r}
\in
C^{0}(\calF^{r})$ such that \[ \beta_{2r} (Y_{1}) = \left\{
\begin{array}{ll} (-1)^{\frac{r(r+1)}{2}} \om_{\la}(X_{1})
\om_{\la}(X_{2}) \cdots \om_{\la}(X_{r}) & \mbox{\rm if~} (Y_{1}) =
(X_{1}) \\ 0 & \mbox{\rm otherwise}. \end{array} \right. \] Thus
$\beta_{2r}=(-1)^{\frac{r(r+1)}{2}}\upsilon(\xi^{*}) \in A^{r}.$ By
construction, $\upsilon$ induces \[
\phi^{-1}: \tilde{H}^{r-1} (F, \C) \longrightarrow H^{r}(A^{\cdot},
\om_{\la}\wedge) \] as long as $\phi$ is an isomorphism.
\end{proof}

For $B\in \nbc(\A)$, define $$ \Xi(B) := \upsilon(\xi(B)^{*}). $$ From
Theorem \ref{fbasis} and Proposition \ref{phiinverse}, we have

\begin{theorem}
\label{XiA}
The set \[ \{\left[\Xi(B)\right]\in
H^{r}(A^{\cdot}, \om_{\la}\wedge) | B\in \bnbc(\A)   \} \] is a basis
for $ H^{r}(A^{\cdot}, \om_{\la}\wedge). $
\end{theorem}

\noindent
Define $\om_{B} := \om_{i_{1} }\cdots \om_{i_{r}}$
for any ordered
base $B = (H_{i_{1} }, \ldots, H_{i_{r}}).$

\begin{remark} \label{ScVsetup1} {In \cite[sections 2, 3]{ScV1}
Schechtman and Varchenko introduced
and studied the flag complex $\{{\mathcal F}^{p}\}$
(not to be confused with the sheaf ${\mathcal F}^{p}$)
and homomorphisms
$S^{p} : {\mathcal F}^{p} \rightarrow A^{p}$.
The homomorphisms $S^{p}$ yield the quasiclassical contravariant
forms.
The abelian group
${\mathcal F}^{r}$ is naturally
identified as a certain factor group of
the group $C_{r-1}(F, \C)$ of $(r-1)$--chains of the Folkman complex
$F$.  Define an epimorphism
$\pi : C^{r-1}(F, \C) \rightarrow {\mathcal F}^{r}$
by $\pi(\xi^{*}) = [\xi]$ for any $\xi \in F^{(r-1)}.$
Then $S^{r}$ and  $\upsilon$
are related by
$\upsilon = S^{r} \circ \pi.$
It is shown in \cite{BrV1} that
set
$\{[\xi(B)] \in {\mathcal F}^{r} | B\in\nbc\}$ is a basis for
${\mathcal F}^{r}$.
Since the homomorphism $S^{r} $ is an isomorphism
by Theorem 3.7 of \cite{ScV1} (recall that we are assuming that
$\la(X)\neq 0$ for all dense $X\in L$),
the set
\begin{eqnarray*}
\{S^{r}([\xi(B)]) | B\in\nbc\}&=&
\{S^{r}\circ \pi(\xi(B)^{*}) | B\in\nbc\}\\&=&
\{\upsilon(\xi(B)^{*}) | B\in\nbc\}\\
&=&
\{\Xi(B) | B\in\nbc\}
\end{eqnarray*}
is a basis for $A^{r}$.

This setup,
together with \cite[Lemma 3.2.5]{ScV1} and Theorem \ref{fbasis},
also yields a different proof of Theorem \ref{XiA} as follows.
It is easy to observe that the flag complex cohomology
 $H^{r} ({\mathcal F}^{\cdot})$ is naturally isomorphic to the
$(r-1)$--dimensional reduced cohomology $\tilde{H}^{r-1}(F) $
of the Folkman complex $F$.
So, the $\bnbc$--bases in Theorem \ref{fbasis} for
$\tilde{H}^{r-1}(F) $
provide a basis for
 $H^{r} ({\mathcal F}^{\cdot}).$
Since the quasiclassical bilinear forms $S^{\bullet}$
give a cochain
isomorphism by Lemma 3.2.5 in \cite{ScV1},
we obtain a basis for $H^{r} (A^{\cdot}, \om_{\la} \we),$ the same
basis as in Theorem \ref{XiA}.
}
\end{remark}

The next theorem requires further restrictions on $\la$.
Let ${\bf P}^{\ell}$ be the complex projective space, which is a
compactification of $V = \C^{\ell}$. Consider the arrangement
$\A_{\infty}$ of projective hyperplanes defined by \[ \A_{\infty} := \{
\overline{H_{1}}, \overline{H_{2}}, \ldots, \overline{H_{\infty}} \},
\]
where $ \overline{H_{i}} $ is the projective closure of $H_{i}~(1\leq
i\leq n)$ and $ \overline{H_{\infty}} := {\bf P}^{\ell} - \C^{\ell}. $
Let $L(\A_{\infty} )$ be the collection of nonempty intersections of
projective hyperplanes in $\A_{\infty}$. Cover ${\bf P}^{\ell}$ by the
standard affine opens $U_{0}, U_{1},\ldots, U_{\ell},$ each of which is
isomorphic to $\C^{\ell}$. Let $\A_{i} ~(0\leq i\leq \ell)$ be the
arrangement in $U_{i} \simeq \C^{\ell}$ obtained by restricting each
projective hyperplane in $\A_{\infty} $ to $U_{i}$. Let $X\in
L(\A_{\infty}) - \{{\bf P}^{\ell} \}$. We say that $X$ is {\em dense}
if
$X\cap U_{i} $ is dense in $\A_{i} $ for $0\leq i \leq \ell$ with
$X\cap
U_{i} \neq \emptyset$. (See \cite{STV1}.) Define $$\la_{\infty} :=
-\sum_{i=1}^{n}  \la_{i}. $$ For $X\in L(\A_{\infty} ) - \{{\bf
P}^{\ell} \}$, let $\la(X)$ be the sum of $\la_{i} $, $i\in\{1, \ldots,
n, \infty\}$, with $X\subseteq \overline{H_{i}}$.

Recall the definitions of
$M$ and the
local system $\calL_{\la}$ from section 1. Combining Theorem
\ref{XiA} with \cite[Theorem 4.1]{STV1} \cite{ESV1}, we obtain

\begin{theorem}
\label{main}
Suppose that none of the $\la(X)$ is a
nonnegative integer for dense $X \in
L(\A_{\infty}) - \{{\bf P}^{\ell}\}$.
Then the set \[ \{\left[\Xi(B)\right]\in H^{r}(M, \calL_{\la}) |
B\in \bnbc(\A)   \} \] is a basis for the local system cohomology
$H^{r}(M, \calL_{\la}). $
\end{theorem}

\begin{example} {Let $\A$ be the arrangement in Example
\ref{example} with defining equation $$Q=(x+1)(x-1)(y+1)(y-1)(x-y),$$
with the linear order defined by the order of the factors in $Q$. Then
\[ \bnbc(\A) = \{ B_{1}, B_{2} \}, \] where $B_{1} = (H_{2}, H_{4})$
and
$B_{2} = (H_{2}, H_{5}).$ We have $\xi(B_{1} ) = ( H_{245} > H_{4} )$
and $\xi(B_{2} ) = ( H_{245}  > H_{5} )$. Thus \[ \Xi(B_{1} ) =
(\la_{2}\om_{2} + \la_{4}\om_{4} + \la_{5}\om_{5}) \la_{4} \om_{4} =
\la_{2}\la_{4}\om_{24} - \la_{4}\la_{5}\om_{45}, \] \[ \Xi(B_{2} ) =
(\la_{2}\om_{2} + \la_{4}\om_{4} + \la_{5}\om_{5}) \la_{5} \om_{5} =
\la_{2}\la_{5}\om_{25} + \la_{4}\la_{5}\om_{45}. \] Let $\la_{\infty} =
- \la_{1} - \cdots -\la_{5}$. Suppose that none of $ \la_{1}, \la_{2},
\la_{3}, \la_{4}, \la_{5}, \la_{2} + \la_{4} + \la_{5}, \la_{1} +
\la_{3} + \la_{5}, \la_{1} + \la_{2} + \la_{\infty}, \la_{3} + \la_{4}
+
\la_{\infty} $ is a nonnegative integer.   Then $\Xi(B_{1} )$ and
$\Xi(B_{2} )$ give a basis for $H^{2}(M, \calL_{\la}).  $}
\end{example}

A different linear order on $\A$ may give another basis for $H^{r}(M,
\calL_{\la}).$ However we have the following. (Compare with the
conjecture in
\cite{Dou1}.)

\begin{proposition} \label{constant} Under the same hypotheses as
Theorem \ref{main}, the transition matrix between two bases for
$H^{r}(M, \calL_{\la})$ obtained from different linear orders on
$\A$ is an integral unimodular matrix independent of $\la$.
\end{proposition}

\begin{proof} The transition matrix is also the transition matrix
between the two bases for $\tilde H^{r-1}(F, \ints). $
\end{proof}

\section{Special cases}

In this section we assume that none of the
$\la(X)$ is a nonnegative integer for dense $X \in
L(\A_{\infty}) -
\{{\bf P}^{\ell} \}$.

Recall \[ \om_{B} = \om_{i_{1} }\cdots \om_{i_{r}} \] for any ordered
base $B = (H_{i_{1} }, \ldots, H_{i_{r}})$. Since $\om_{B}$ is a
``monomial'' while $\Xi(B)$ is a sum of several $\om_{B}$, it
might be desirable that the set $\{\left[\om_{B}\right] ~|~ B\in
\bnbc(\A)\}$ form a basis for $H^{r}(M, \calL_{\la})$. This will not
hold in general, even for arrangements of rank two.
In this section,
we will prove that these monomial forms give a basis in two special
situations:

\begin{itemize}
\item when the linear order on $\A$ is unmixed (Def. \ref{unmixeddef}),
\item when the linear order on $\A$ is admissible (Def.
\ref{admissibledef}) and $r=2.$
\end{itemize}

\noindent
The first case contains arrangements in general position,
arrangements in general position to infinity, and normal arrangements
among others.

\begin{definition}
\label{unmixeddef}
{For each maximal $X\in L$, let \[ \nbc_{X} := \{B \in\nbc ~|~
\cap B = X\}. \] We say that $X$ is {\em unmixed} if either
$\nbc_{X}\subseteq \bnbc$ or $\nbc_{X}\cap \bnbc = \emptyset.$
A linear order on $\A$ is {\em unmixed} if any maximal $X\in L$ is
unmixed.}
\end{definition}

\begin{theorem}
\label{unmixedtheorem}
If the linear order on $\A$ is
unmixed, then the set $\{\left[\om_{B}\right] ~|~ B\in \bnbc(\A)\}$
is a basis for $H^{r}(M, \calL_{\la})$.
\end{theorem}

\begin{proof} Let $X\in L$ be an arbitrary maximal element with
$\nbc_{X} \subseteq \bnbc.$  Then each $\Xi(B), ~B\in \nbc_{X},$
belongs to
$$ \sum_{B'\in \nbc_{X}} \C \om_{B'} \subseteq \sum_{B'\in \bnbc} \C
\om_{B'}.$$ The result follows from Theorem \ref{main} and a dimension
argument.
\end{proof}

\begin{example}
\label{}
{An arrangement $\A$ is said to be {\em in
general position} if (i) $n = |\A| \geq \ell + 1$, (ii) $\codim
(H_{i_{1}}\cap \cdots \cap H_{i_{k}}) = k$ whenever $1\leq k\leq \ell,$
and (iii)    $H_{i_{1}}\cap \cdots \cap H_{i_{k}} = \emptyset$ whenever
$k
> \ell$. Let $X\in L$ be a maximal element.  Then $\nbc_{X} $ is a
singleton.  Thus any linear order on $\A$ is unmixed. We have \[
\bnbc(\A) = \{(H_{i_{1}}, \ldots, H_{i_{\ell}}) ~|~ 1 < i_{1} < \cdots
<
i_{\ell} \leq n \}, ~~\mbox{and}\] \[ \Xi(B) = \la_{i_{1}}\cdots
\la_{i_{\ell}}
\om_{B}, \] where $B=(H_{i_{1}}, \ldots, H_{i_{\ell}}) \in\bnbc(\A)$.
This basis
coincides with the basis constructed in \cite[p.292]{Aom1}.}
\end{example}

\begin{example} \label{} {An arrangement $\A$ is said to be {\em
normal} \cite[1.4]{Var1} if $|\A_{X}| = \codim(X)$ for all $X\in
L(\A)$.
(It is equivalent to say that $\bigcup \A$ is a normal crossing divisor
in $V$.) Arrangements in general position are normal. Let $X\in L$ be a
maximal element.  Then $\nbc_{X} $ is again a singleton.  Thus any
linear order on $\A$ is unmixed. We have \[ \Xi(B) = \la_{i_{1}}\cdots
\la_{i_{r}} \om_{B}, \] where $B=(H_{i_{1}}, \ldots, H_{i_{r}})
\in\bnbc(\A)$.}
\end{example}

\begin{example} \label{} {An arrangement $\A$ is said to be {\em in
general position to infinity} \cite[6.2]{Yuz1} if $H_{i_{1}}\cap \cdots
\cap H_{i_{k}} \neq \emptyset$ whenever $1 \leq k \leq \ell$.
This implies that there are no parallels among the elements of $L$, so
that
the hyperplane at infinity is
generic relative to $\A$. In particular, general position arrangements
are in
general position to infinity.

Let $X\in L$ be a
maximal element.  Note that $\A_{X}$ is a central generic arrangement
\cite[5.22]{OrT1}.
We have $\nbc_{X}\cap \bnbc = \emptyset$ if and only if
$X\subseteq H_{1}$. Otherwise  $\nbc_{X}\subseteq \bnbc.$ Thus any
linear order on $\A$ is unmixed. We have \[ \bnbc(\A) = \{(H_{i_{1}},
\ldots, H_{i_{r}}) ~|~ 1 < i_{1} < \cdots < i_{r} \leq n,
i_{1} = \min \A_{X},
X = H_{i_{1}\cdots i_{r}}\}.
\] }
\end{example}

\begin{example} \label{} {Let $H_{1} \in \A$ be generic,
that is, $H_{1} $ transversely intersects $Y$ unless $Y\in
L(\A - \{H_{1} \})$ is maximal in $L(\A - \{H_{1} \})$. Let $X\in L$ be
a maximal element.  Then $\nbc_{X}\cap \bnbc = \emptyset$ if and only
if
$X\subseteq H_{1}$. Otherwise  $\nbc_{X}\subseteq \bnbc.$ Thus any
linear order on $\A$ in which $H_{1}$ is the first hyperplane is
unmixed. In this case, we have \[ \bnbc(\A) = \{(H_{i_{1}}, \ldots,
H_{i_{r}})\in \nbc(\A) ~|~ 1 < i_{1} < \cdots < i_{r} \leq n \}. \] }
\end{example}

\begin{remark} {Suppose that $\A$ is complexified real and $r=\ell$.
Then the number of bounded chambers of the real form of $\A$ is equal
to
$|\bnbc|$.  In this case Varchenko \cite[6.2]{Var1} associated a
differential $\ell$--form $\eta_{\Delta}\in A^{\ell}(\A)$ (which is not
necessarily a monomial) to each bounded chamber $\Delta$. Recall
the definition of
$U_{\la}$ from section 1 and define $\om_{\Delta} = U_{\la}
\eta_{\Delta}.$ The form $\om_{\Delta} $ is called the hypergeometric
form associated to  $\Delta.$ The hypergeometric integrals give a
determinant $\det \left[ \int_{\Gamma} \om_{\Delta}\right].$ The main
results in \cite{Var1} are beautiful formulas for the determinant when
$\A$ is in general position (Thm.1.1), normal (Thm.1.4), or in general
position to infinity (Thm.6.1).  In particular, these determinants
are nonzero, so that in these cases the set $\{\eta_{\Delta}\}$ gives a
basis
for $H^{\ell}(M, \calL_{\la}).$ The relationship between
$\{\eta_{\Delta}\}$ and $\{\Xi(B)\}$ is intriguing. They are not the
same in general but coincide when $\A$ is normal or in general position
to infinity.} \label{Varch} \end{remark}

\noindent
Finally we specialize to the case $r=2$.
\begin{definition}
\label{admissibledef} {The linear order on $\A$ is called {\em
admissible} if there exists an integer $\nu$ such that $H_{i}$ and
$H_{1} $ are parallel if and only if $1\leq i < \nu.$}
\end{definition}

If the linear order on $\A$ is admissible, then it is not difficult to
see \[ \bnbc (\A) = \{ (H_{i_{1} }, H_{i_{2} }) \in \nbc ~|~ 1 < i_{1}
<
i_{2} \neq \nu\}. \]

\begin{proposition} \label{admissibletheorem} Suppose $\A$ is an
arrangement
of rank 2 with admissible linear order.  Then the set
$\{ \left[\om_{B} \right]~|~B\in \bnbc\}$
gives a basis for $H^2 (M, \calL_{\la}).$
\end{proposition}

\begin{proof} Since $\dim H^{2}(M, \calL_{\la}) = |\bnbc|,$ it suffices
to show that the set $\{ \left[\om_{B} \right] ~|~B\in\bnbc\}$ spans
$H^{2}(M, \calL_{\la})$. Define \[ N := \sum_{B\in\bnbc} \C \om_{B} +
d_{\la} (A^{1}). \] We want to show $N = A^{2}. $ By Theorem
\ref{main},
it is enough to show that $\Xi(B) \in N$ for all $B\in\bnbc$. Let
$B\in\bnbc$ and $X = \cap B$. If $X$ is unmixed, then \[ \Xi(B) \in
\sum_{B'\in\nbc_{X}} \C \om_{B'} \subseteq \sum_{B'\in\bnbc} \C
\om_{B'} \subseteq N. \] Suppose that $X$ is mixed.  Then $B = (H_{i},
H_{j}) \in \bnbc$ with $1 < i < \nu < j$ and $X = H_{i} \cap H_{\nu} =
H_{i} \cap H_{j}$. Note that $(H_{i}, H_{p}) \in\bnbc$ for all $p >
\nu.$ Thus $\omega_{p}\omega_{q} = \omega_{i}\omega_{q} -
\omega_{i}\omega_{p} \in N   $ if $\nu \not\in \{p, q\}$ and $H_{p}
\cap
H_{q} = X$. Therefore we have the following congruence relations modulo
$N$: \begin{eqnarray*} \Xi(B) &=& \om_{\la}(X) \wedge \la_{j}\om_{j}
\equiv \la_{\nu}\om_{\nu} \wedge \la_{j} \om_{j}\\ &=& \la_{\nu}\la_{j}
(\om_{ij} - \om_{i \nu}) \equiv - \la_{\nu}\la_{j} \om_{i \nu} \equiv
\la_{j} d_{\la}(\om_{i}) \equiv 0. \end{eqnarray*} This proves
$\Xi(B)\in N.$ \end{proof}

\noindent
Proposition \ref{admissibletheorem} was independently proved by
M. Kita \cite{Kit1}.  He gives
 a direct proof
which doesn't use Theorem \ref{main}.

If the linear order is not admissible, the
set $\{ \left[\om_{B} \right] ~|~B\in\bnbc\}$ does not give a basis for
$H^{2} (M, \calL_{\la})$ in general unless we impose additional
unnatural genericity conditions on $\la$.

\end{section}

\begin{ack}
We thank M.~Kita, P.~Orlik, A.~Varchenko and G.~Ziegler
for helpful discussions.    The second author is also grateful to
H.~Kanarek
for a helpful and inspiring conversation describing his work.
\end{ack}


\begin{thebibliography}{99}

\bibitem{Aom1} Aomoto, K.: Les \'equations aux diff\'erences
lin\'eaires
et les int\'egrales des fonctions multiformes, J. Fac. Sci., Univ.
Tokyo, {\bf 22} (1975) 271--297

\bibitem{Aom2} Aomoto, K.: On vanishing of cohomology attached to
certain many valued meromorphic functions, J. Math. Soc. Japan, {\bf
27}
(1975), 248--255

\bibitem{AoK1} Aomoto, K., Kita, M.: Hypergemetric functions (in
Japanese), Springer-Verlag, Tokyo, 1994

\bibitem{B1}
A.~Bj\"orner:
On the homology of geometric lattices,
Algebra Universalis {\bf 14} (1982), 107-128.

\bibitem{B2}
A.~Bj\"orner:
Homology and shellability of matroids and geometric lattices,
in Matroid Applications (ed.\ N.~White),
Cambridge University Press 1992, pp.~226-283.

\bibitem{BjZ1} Bj\"orner, A., Ziegler, G.: Broken circuit complexes:
Factorizations and generalizations, J. Combin. Theory Ser. B {\bf 51}
(1991), 96-126.

\bibitem{BrV1} Brylawski, T., Varchenko, A.: The determinant formula
for
a matroid bilinear form, preprint, 1994.

\bibitem{Cra} Crapo, H.:  A higher invariants for matroids, J. of
Combinatorial Theory {\bf 2} (1967) 406--417.

\bibitem{Dou1} Douai, A.: Autour d'une conjecture de Varchenko sur les
int\'egrales hyperg\'eom\'etriques: le point de vue des syst\`emes
d'\'equations aux diff\'erences finies, preprint, 1994.

\bibitem{ESV1} Esnault, H., Schechtman, V.,  Viehweg, E.: Cohomology of
local systems of  the complement of hyperplanes, Invent. math.  {\bf
109} (1992), 557-561; Erratum, ibid.  {\bf 112} (1993) 447.


\bibitem{Gel1} Gelfand, I.M.: General theory  of hypergeometric
functions, Soviet Math. Dokl.  {\bf 33} (1986), 573--577.

\bibitem{God1} Godement,  R.: Topologie alg\'ebrique et th\'eorie des
faisceaux,  Hermann, Paris, 1958.

\bibitem{Kan1} Kanarek, H.: Ph.D. thesis, University of Essen, in
preparation.

\bibitem{Kit1} Kita, M.: Private communication.

\bibitem{Koh1} Kohno, T.: Homology of a local system on the complement
of hyperplanes, Proc. Japan Acad. {\bf 62}, Ser. A (1986), 144--147.


\bibitem{OrT1} Orlik, P., Terao, H.: Arrangements of  Hyperplanes.
Grundlehren der Math. Wiss. {\bf 300},   Springer Verlag, 1992.


\bibitem{STV1} Schechtman, V. V., Terao, H., Varchenko, A. N.:
Cohomology of local systems and the Kac-Kazhdan condition for singular
vectors, preprint, 1994.

\bibitem{ScV1} Schechtman, V. V., Varchenko, A. N.: Arrangements of
hyperplanes and Lie algebra homology, Invent. math. {\bf 106} (1991),
139--194.

\bibitem{Var1} Varchenko, A.N.: The Euler beta-function,  the
Vandermonde determinant, Legendre's equation,  and critical values of
linear functions on a  configuration of hyperplanes. I. Math. Ussr
Izvestiya,  {\bf 35} (1990), 543-571, II. Math. Ussr Izvestiya, {\bf
36}
(1991), 155--167.

\bibitem{Var2} Varchenko, A.N.: Multidimensional hypergeometric
functions
and  representation theory of  Lie algebras
and quantum groups, Advanced Series in Mathematical Physics -
Vol. 21,  World Scientific Publishers, to appear.

\bibitem{WaW} Wachs, M. L., Walker, J.W.:
On geometric semilattices, Order, {\bf 2} (1986), 367--385.

\bibitem{Yuz1} Yuzvinsky, S.: Cohomology of the Brieskorn-Orlik-Solomon
algebras, preprint, 1994.

\bibitem{Zie1} Ziegler, G.: Matroid shellability, $\beta$--systems, and
affine arrangements, J. of Alg. Combinatorics, {\bf 1} (1992),
283--300.



\end{thebibliography}
\end{document}